\title{Regularity Criteria for Navier-Stokes Equations with Slip Boundary
Conditions on Non-flat Boundaries via Two Velocity Components}
\author{Hugo Beir\~{a}o da Veiga$^{1,}$ \footnote{Partially supported  by FCT (Portugal) under grant UID/MAT/04561/3013.}\qquad Jiaqi Yang$^{2,}$\footnote{Hugo Beir\~{a}o da Veiga (\texttt{bveiga@dma.unipi.it}) and Jiaqi Yang (\texttt{yjq@imech.ac.cn})}}
\date{
\small $^1$ Department of Mathematics, Pisa University, Pisa, Italy\\
\small $^2$ Key Laboratory for Mechanics in Fluid Solid Coupling Systems, Institute of Mechanics, Chinese Academy of Sciences, Beijing 100190, China}
\documentclass[11pt]{article}
\usepackage{amsfonts}
\usepackage{mathrsfs}
\usepackage{amsmath}
\usepackage{amssymb,extarrows}
\usepackage{multicol}
\usepackage{float}
\usepackage{makeidx}
\usepackage{layout}
\usepackage{array}
\usepackage{a4wide}
\usepackage{boxedminipage}
\usepackage{hyperref}
\usepackage{latexsym}
\usepackage{color}
\usepackage{dsfont}
\usepackage{graphicx}
\usepackage{ulem}
\usepackage{amsthm}

\newtheorem{theorem}{Theorem}[section]

\newtheorem{assumption}[theorem]{Assumption}
\newtheorem{problem}[theorem]{Problem}
\theoremstyle{remark}
\newtheorem{remark}{Remark}[section]

\theoremstyle{definition}

\numberwithin{equation}{section}

\newcommand{\p}{\partial}
\newcommand{\e}{\epsilon}

\newcommand{\R}{\mathbb{R}}

\newcommand{\f}{\frac}
\newcommand{\n}{\nabla}

\newcommand{\bom}{\boldsymbol{\omega}}
\newcommand{\bu}{\mathbf{u}}
\newcommand{\bn}{\mathbf{n}}
\newcommand{\ed}{\end{document}}
\begin{document}
\maketitle

\begin{abstract}
H.-O. Bae and H.J. Choe, in a 1997 paper, established a regularity criteria for the incompressible Navier-Stokes equations in the whole space $\R^3$ based on two velocity components. Recently, one of the present authors extended this result to the half-space case $\R^3_+\,.$ Further, this author in collaboration with J. Bemelmans and J. Brand extended the result to cylindrical domains under physical slip boundary conditions. In this note we obtain a similar result in the case of smooth arbitrary boundaries, but under a distinct, apparently very similar, slip boundary condition. They coincide just on flat portions of the boundary. Otherwise, a reciprocal reduction between the two results looks not obvious, as shown in the last section below.
\end{abstract}

\noindent \textbf{Mathematics Subject Classification:} 35Q30.

\vspace{0.2cm}
\noindent \textbf{Keywords:} Navier-Stokes equations; Slip boundary conditions; No flat boundaries; Two components regularity criterium.

\vspace{0.2cm}

\section{Introduction}
The starting point of the present paper is the well known Prodi-Serrin (P-S) sufficient condition for regularity of the solutions to the incompressible Navier-Stokes equations
\begin{equation}\label{eq:NS}
\begin{cases}
\p_t\bu+\bu\cdot\n\bu-\Delta\bu+\n p=0\,,\quad&\text{in $\Omega\times(0,T]$}\,,\\
\n\cdot\bu=0\,,\quad&\text{in $\Omega\times(0,T]$}\,.
\end{cases}
\end{equation}
where $\bu=(u_1,u_2,u_3)$ denotes the unknown velocity of the fluid and $p$ the pressure.
To immediately set limits to the circle of our interests, assume for now on that $\Omega \subset \,\R^3\,$ is a bounded, smooth domain, even if many results quoted below hold for larger space dimensions. For the time being, assume that suitable boundary conditions are imposed to the velocity $\bu$ .\par%
The global existence of the so called \emph{weak} solutions to system \eqref{eq:NS} goes back to J. Leray \cite{Leray} and E. Hopf \cite{Hopf} classical references. See also A.A. Kiselev and O.A. Ladyzhenskaya \cite{kiselev}, and J.L. Lions \cite{lions}. Below, solutions of \eqref{eq:NS} are intended in this sense.%

\vspace{0.2cm}

A main classical open mathematical problem is to prove, or disprove, that weak solutions are necessarily strong under reasonable but general assumptions, where strong means that
\begin{equation}\label{eq:Strong}
\bu\in L^{\infty}(0,T;H^1(\Omega))\cap L^2(0,T;H^2(\Omega))\,.
\end{equation}
In this context, a remarkable and classical sufficient condition for uniqueness and regularity is the so-called Prodi-Serrin condition, P-S in the sequel, namely
\begin{equation}\label{eq:P-S}
\bu\in L^{q}(0,T;L^p(\Omega))\,,\quad \f{2}{q}+\f{3}{p}=\,1\,,\quad p>3\,.
\end{equation}

\vspace{0.2cm}

Concerning condition \eqref{eq:P-S}, we transcribe from \cite{Beirao18}, Section 1, the following considerations: Assumption \eqref{eq:P-S} was firstly considered by
G.~Prodi in his paper \cite{Prodi} of 1959.
He proved uniqueness under this last assumption. See also C.~Foias, \cite{Foias}.
Furthermore, J.~Serrin, see \cite{Serrin62} and \cite{Serrin63},
particularly proved interior spatial regularity under the stronger (non-strict) assumption
\begin{equation}\label{eq:lps2}
\bu\in L^q(0,T; L^p(\Omega)),\quad\frac2q + \frac 3p<1,\quad p>3.
\end{equation}
Concerning the above problems, see also O.A.~Ladyzhenskaya's contributions
\cite{lady-2} and \cite{ladyz}.
The above setup led to the nomenclature Prodi-Serrin condition.\par%
Complete proofs of the strict regularity result (i.e. under assumption \eqref{eq:P-S})
were given by H.~Sohr in \cite{sohr}, W.~von~Wahl in \cite{vonwahl}, and Y.~Giga in \cite{giga}.
A simplified version of the proof was given in reference \cite{MR91e:35066}, to which we refer also for bibliography.
For a quite complete overview on the main points, and references, on the initial-boundary value problem for Navier-Stokes equations we strongly recommend Galdi's contribution \cite{galdi2}. Further, we refer to \cite{prodi-2} and \cite{Serrin63}, as sources for information on the historical context of the P-S condition by the initiators themselves. \par%
Finally, we recall that L.~Escauriaza, G.~Seregin, and V.~\v{S}ver\'ak, see \cite{seregin},
extended the regularity result to the case $ (q,p)=(\infty,3).$\par%

\vspace{0.2cm}

A significant improvement of the P-S condition was obtained by H.-O. Bae and H.J. Choe \cite{Bae97}, see also \cite{Bae07}. They proved, in the whole space case, that it is sufficient for regularity of solutions that two components of the velocity satisfy the above condition \eqref{eq:P-S}. For convenience we call here this situation as being the \emph{restricted}
P-S condition. In 2017, one of the authors, see \cite{Beirao17}, extended this result to the half-space $\R^3_+$ under slip boundary conditions. In this case, the truncated 2-dimensional vector field $\bar{\bu}$ cannot be chosen arbitrarily. The omitted component has to be the normal to the boundary.\par%
Very recently, in reference \cite{Beirao18}, the result was extended to a cylindrical type three-dimensional domain, consisting on the complement set between two co-axial circular cylinders, with radius $\rho_0$ and $\rho_1\,,$  $\,0<\,\rho_0 <\,\rho_1,$ periodic in the axial direction, under the physical slip boundary condition
\begin{equation}\label{eq:bbb}
\bu\cdot\bn=0\,, \quad [D(\bu)\bn]\cdot\tau=\,0\,,\quad\text{on $\p\Omega$}\,,
\end{equation}
where $D(\bu)=\f{\n\bu+(\n\bu)^T}{2}$ is the shear stress. The above exclusion of an interior cylinder was done to avoid the radial coordinate singularities on the symmetry axis, which consideration is out of interest in our context. Below we obtain a similar result, extended to domains with general non-flat boundaries, but under the slip boundary condition \eqref{eq:bry}. The two boundary conditions coincide just on flat portions of the boundary. Otherwise, a reciprocal reduction between the two results looks not obvious. This claim is shown in the last section.

\vspace{0.2cm}

Again by following \cite{Beirao18} we recall that after the contribution by H.-O. Bae and H.J. Choe, related papers appeared that
particularly concerned assumptions on two components of velocity or vorticity,
see \cite{wolf}, \cite{bdvjmfm}, \cite{Beirao17}, \cite{berselli}, and \cite{chae-choe}.
There are also many papers dedicated to sufficient conditions for regularity
which depend merely on one component, see, for instance,
\cite{cao-titi},  \cite{cheng}, \cite{kukavica},  \cite{neustupa}, \cite{zhang},
and \cite{zhou-pokorny}.

\vspace{0.3cm}

Before going on we want to motivate the particular choice of the domain made below. It takes into account that the real significance of the result has essentially a local character. First of all, a global regular (i.e., without singularities) system of coordinates, two of them parallel and the third orthogonal to the boundary, does not exist in general, even in an arbitrarily thin neighbourhood of the full boundary, as in the case of a sphere and even in the case of a spherical corona. In fact, singularities typically appear, like on the above two cases, and even in full cylinders (due to the symmetry axis). The cylindrical case considered in reference \cite{Beirao18} is an exception (see below) due to the removal of a neighbourhood of the symmetry axis.\par%
Luckily, the above type of coordinates' system exists in sufficient small neighbourhoods of any regular boundary point. Hence, to illustrate the full significance of our thesis in a simple, but still convincing way, it looks sufficient to prove it near any ``small'' piece of smooth boundary with an \emph{arbitrary geometrical shape}. This is our aim below. The restrictions on the domain $\Omega$ below are made in accordance with these lines, a choice which covers the very basic situation, in the simplest way.%

\section{Main Results}
In the sequel we assume the slip boundary condition
\begin{equation}\label{eq:bry}
\bu\cdot\bn=0\,,\quad \bom\times\bn=0\,\quad\text{on $\p\Omega$}\,,
\end{equation}
where $\bom=\n\times\bu$ is the vorticity, $\bn$ is the outward normal of $\p\Omega$, and $\Omega\in\R^3$ is a smooth domain satisfying the following condition:
\begin{assumption}\label{assump}
\textit{
There exists a curvilinear orthogonal system of coordinates
\[
q(x)=(q_1(x),q_2(x),q_3(x))
\]
such that $\Omega$ can be transited into
\[
\hat{\Omega}\triangleq\{(q_1,q_2,q_3): 0\leq q_1<1,\, 0\leq q_2<1,\, 0<\rho_0\leq q_3\leq\rho_1\}\,,
\]
where the axis $q_3$ direct to the outward normal on the boundary $\p\hat{\Omega}^1:=\{(q_1,q_2,q_3): q_3=\rho_1\}$ (the inward normal on the boundary $\p\hat{\Omega}^0:=\{(q_1,q_2,q_3): q_3=\rho_0\}$, respectively),  and $q_1$, $q_2$ are periodic. 
}
\end{assumption}
\begin{remark}
The above "small" piece of a generical smooth boundary is here represented by $q_3=\,\rho_0\,,$ and  $q_3=\,\rho_1\,.$\par%
\end{remark}
\begin{remark}\label{re:due}
It is worth noting that the slip boundary condition  \eqref{eq:bry} is equivalent to
\begin{equation}\label{eq:Navier-bry-ka}
\bu\cdot\bn=0\,,\quad [D(\bu)\bn]\cdot\tau=\,-\kappa_{\tau}\, \bu\cdot\tau\,,
\end{equation}
where $\tau$ stands for any arbitrary unit tangential vector on $\p\Omega\,,$ and $\kappa_{\tau}$ is the principal curvature in the $\tau$ direction, positive if the center of curvature lies inside $\Omega\,.$\par%
The above claim follows immediately by appealing  to equation  (5.2) in \cite{beirao-crispo}, namely
\begin{equation}\label{eq:BV-CR}
\quad [D(\bu)\bn]\cdot\tau=\,\frac{1}{2}\,(\bom\times\bn)\cdot\tau -\,\kappa_{\tau}\,\bu\cdot\tau\,.
\end{equation}
For a mathematical treatment of some aspects related to slip boundary conditions imposed on smooth, but generic, boundaries see also \cite{bdvadvances}, and the pioneering paper \cite{so-sca}.
\end{remark}

\vspace{0.2cm}

Next we recall some facts on curvilinear coordinates. The Lam\'e coefficients (scale factors) of the transition to the system of coordinates $q$ are denoted by the letters $H_i$
\begin{equation*}
H_i(q)=\left(\sum_{j=1}^3\left(\f{\p x_j}{\p q_i}\right)^2\right)^{\f12}\,,\quad \text{$i=1\,,2\,,3$.}
\end{equation*}
Let $\hat{e}_i=\f{1}{H_i}\f{\p x}{\p q_i}$, $i=1\,,2\,,3$. Note that $\,|\hat{e}_i|=\,1$ and $\,\hat{e}_i\cdot\n=\f{1}{H_i}\f{\p}{\p q_i}\,.$  One can write
\[
\bu(x)=\hat{\bu}(q)=\hat{u}_1(q)\hat{e}_1+\hat{u}_2(q)\hat{e}_2+\hat{u}_3(q)\hat{e}_3
\]
and
\[
\bom(x)=\hat{\bom}(q)=\hat{\omega}_1(q)\hat{e}_1+\hat{\omega}_2(q)\hat{e}_2+\hat{\omega}_3(q)\hat{e}_3.
\]
It is well known (see for example \cite{KIR} and \cite{batchelor}) that
\begin{equation}\label{eq:div}
\n\cdot\bu=\f{1}{H_1H_2H_3}\left(\f{\p(\hat{u}_1H_2H_3)}{\p q_1}+\f{\p(\hat{u}_2H_1H_3)}{\p q_2}+\f{\p(\hat{u}_3H_1H_2)}{\p q_3}\right)
\end{equation}
and
\begin{equation}\label{eq:curl}
\begin{split}
\n\times\bu=&\f{1}{H_2H_3}\left(\f{\p(\hat{u}_3H_3)}{\p q_2}-\f{\p(\hat{u}_2H_2)}{\p q_3}\right)\hat{e}_1\\
&+\f{1}{H_1H_3}\left(\f{\p(\hat{u}_1H_1)}{\p q_3}-\f{\p(\hat{u}_3H_3)}{\p q_1}\right)\hat{e}_2\\
&+\f{1}{H_1H_2}\left(\f{\p(\hat{u}_2H_2)}{\p q_1}-\f{\p(\hat{u}_1H_1)}{\p q_2}\right)\hat{e}_3.
\end{split}
\end{equation}
We state our main result as follows.
\begin{theorem}\label{thm}
Let $\Omega$ satisfy Assumption \ref{assump}, and suppose that there exist two positive constants $c$ and $C$ such that
\begin{equation}\label{eq:assumption}
c\leq H_i\leq C\quad \text{and}\quad\left|\f{\p^2 x_i}{\p q_i\p q_j}\right|\,,\left|\f{\p^3 x_i}{\p q_i\p q_j\p q_k}\right|\leq C\,,
\end{equation}
for any $i\,,j\,,k=1\,,2\,,3$.  Let $\bu$ be a weak solution of the system \eqref{eq:NS} under the boundary condition \eqref{eq:bry}, and set $\bar{\bu}=\hat{u}_1\hat{e}_1+\hat{u}_2\hat{e}_2$. If $\bar{\bu}$ satisfies
\begin{equation}\label{eq:RCcondition}
\bar{\bu}\in L^q(0,T;L^p(\Omega))\,,\quad\f{2}{q}+\f{3}{p}\leq1,\quad p>3\,,
\end{equation}
then the solution $\bu$ is strong, namely,
\[
\bu\in L^{\infty}(0,T;H^1(\Omega))\cap L^2(0,T;H^2(\Omega))\,.
\]
\end{theorem}
Note that assumption  \eqref{eq:assumption} implies $\, |\p_iH_j|\,,|\p_{ij}H_k|\leq C\,.$\par%
It is worth noting that our proof applies to a more general set of geometrical situations. let's just give some hint in this direction.
\begin{remark}\label{jordan}
The above statement does not contain the result proved in reference \cite{Beirao18}, due to the distinct boundary conditions, see Section \ref{moregen}. On the other hand, we may replace the two circular, vertical, cylinders by more general vertical cylinders where the external circle $q_3=\,\rho_1\,$ is replaced by a smooth Jordan curve $\,\gamma_1\,$, and the internal circle  $q_3=\,\rho_0$ by a parallel Jordan curve $\,\gamma_0\,$, at a sufficient small distance $\,\delta>0\,$ from $\gamma_1\,$. The coordinate $\theta$ is now an arc length coordinate on $\gamma_1$. All points in the same orthogonal segment to $\gamma_0\,$ and $\gamma_1\,$ enjoy the same $\theta$ coordinate. The coordinate $\,r\in \,(0,\,\delta)\,$ is given by the distance to $\gamma_1$. The ``vertical" coordinate $\,z\,$ preserves his periodic character. Clearly, the role played by the above Jordan curve may be immediately extended to much more general situations.\par%
Another significant application is obtained by replacing the above two cylindrical boundaries by two torus of revolution, generated by revolving two concentric circles $\,\gamma_0\,$ and $\gamma_1\,$ about an axis coplanar with the circles, which does not touch the circles (roughly, we obtain the complement set between two closed tubes). Now $\,z \in \,[0,\,2\,\pi)\,$ is an angular periodic coordinate, the toroidal coordinate. The result still applies by replacing the two circles by two parallel Jordan curves.%
\end{remark}
Let's propose the following benchmark problem:
\begin{problem}\label{problema}
Consider two concentric spheres $\,\Omega_R$ and  $\,\Omega_\rho\,,$ of radius respectively $\rho$ and $R$,  $\,0<\,\rho<\,R\,$. Let $\,u\,$ be a weak solution in $\,\Omega_R\times(0,T]\,$ of \eqref{eq:NS} under one of the above slip boundary conditions. Further, assume that the restricted P-S condition holds in  $\,(\Omega_R-\,\Omega_\rho)\,\times(0,T]\,$  with respect to the tangential components of the velocity, and holds in $\,\Omega_\rho\,\times(0,T]\,$ with respect to two
arbitrary components of the velocity. Problem: To prove that $u$ is a strong solution in $\,\Omega_R\times(0,T]\,$.\par%
\end{problem}
\section{Proof of Theorem \ref{thm}}\label{seq:proof}
\begin{proof}
We start by reducing the system \eqref{eq:NS} under the boundary condition \eqref{eq:bry} into the classical vorticity form
\begin{equation*}
\begin{cases}
\p_t\bom+\bu\cdot\n\bom-\bom\cdot\n\bu-\Delta\bom=0\,,\quad&\text{in $\Omega\times(0,T]$}\,,\\
\n\cdot\bu=0\,,\quad&\text{in $\Omega\times(0,T]$}\,,\\
\bu\cdot\bn=0\,,\quad \bom\times\bn=0\,,\quad&\text{on $\p\Omega$}\,.
\end{cases}
\end{equation*}
Then we take the scalar product with $\bom$, and integrate by parts. One easily gets
\begin{equation}\label{eq:Basicestimate}
\begin{split}
\f12\p_t\int_{\Omega}|\bom|^2dx+\int_{\Omega}|\n\bom|^2dx=\int_{\p\Omega}\bn\cdot\n\bom\cdot\bom dS+\int_{\Omega}\bom\cdot\n\bu\cdot\bom dx
:=I_1+I_2\,.
\end{split}
\end{equation}
Next, we focus on the estimates of $I_1$ and $I_2$.
\vskip 0.1 true cm
\noindent\textit{Control of $I_1$:}
First, it follows from \eqref{eq:bry} that
\begin{equation}\label{eq:bry1}
\hat{u}_3=0\,,\quad \hat{\omega}_1=\hat{\omega}_2=0\,, \quad\text{as \, $q_3=\rho_0\,,\rho_1$}\,.
\end{equation}
Let $\p\Omega^l=\p\hat{\Omega}^l:=\{(q_1,q_2,q_3)\in\hat{\Omega}: q_3=\rho_l\,\}$\,, where $ l=0\,,1$. One can deduce from \eqref{eq:bry1} that
\begin{equation}\label{eq:comes}
\begin{split}
(2l-1)\int_{\p\Omega^l}\bn\cdot\n\bom\cdot\bom dS=&(2l-1)\int_{\p\Omega^l}\bn\cdot\n\left(\f{|\bom|^2}{2}\right) dS\\
=&\int_0^1\int_0^1\left[\p_{q_3}\left(\f{|\hat{\bom}|^2}{2}\right)H_1H_2H_3^{-1}\right]\Big|_{q_3=\rho_l}dq_1dq_2\\
=&\int_0^1\int_0^1\left[(\p_{q_3}\hat{\omega}_3)\,\hat{\omega}_3H_1H_2H_3^{-1}\right]|_{q_3=\rho_l}dq_1dq_2\\
=&\int_0^1\int_0^1\left[\p_{q_3}\left(H_1H_2\hat{\omega}_3\right)H_3^{-1}\hat{\omega}_3\right]|_{q_3=\rho_l}dq_1dq_2\\
&-\int_0^1\int_0^1\left[\p_{q_3}\left(H_1H_2\right)H_3^{-1}\hat{\omega}^2_3\right]|_{q_3=\rho_l}dq_1dq_2\,.
\end{split}
\end{equation}
Since $\n\cdot\bom=0$, from \eqref{eq:div} one gets
\begin{equation*}
\f{\p(\hat{\omega}_1H_2H_3)}{\p q_1}+\f{\p(\hat{\omega}_2H_1H_3)}{\p q_2}+\f{\p(\hat{\omega}_3H_1H_2)}{\p q_3}=0\,,
\end{equation*}
which gives
\begin{equation*}
\begin{split}
&\int_0^1\int_0^1\left[\p_{q_3}\left(H_1H_2\hat{\omega}_3\right)H_3^{-1}\hat{\omega}_3\right]|_{q_3=\rho_l}dq_1dq_2\\
=&-\int_0^1\int_0^1\left[\p_{q_1}\left(H_2H_3\,\hat{\omega}_1\right)H_3^{-1}\hat{\omega}_3\right]|_{q_3=\rho_l}dq_1dq_2
-\int_0^1\int_0^1\left[\p_{q_2}\left(H_1H_3\hat{\omega}_2\right)H_3^{-1}\hat{\omega}_3\right]|_{q_3=\rho_l}dq_1dq_2\\
=&\,0\,,
\end{split}
\end{equation*}
since $\,\hat{\omega}_1=\,\hat{\omega}_2 =\,\p_{q_1}\,\hat{\omega}_1=\,\p_{q_2}\hat{\omega}_2=\,0\,$ on $\,\p\hat{\Omega}^l\,.$ Hence, one obtains
\begin{equation*}
\begin{split}
(2l-1)\int_{\p\Omega^l}\bn\cdot\n\bom\cdot\bom dS
=-\int_0^1\int_0^1\left[\p_{q_3}\left(H_1H_2\right)H_3^{-1}\hat{\omega}^2_3\right]|_{q_3=\rho_l}dq_1dq_2
\end{split}
\end{equation*}
By appealing to \eqref{eq:assumption} one shows that
\[
\left|\int_{\p\Omega}\bn\cdot\n\bom\cdot\bom dS\right|\leq C\int_{\p\Omega}|\bom|^2dS\leq\,\|\,|\bom\,|^2\|_{W^{1,\,1}(\Omega)}\,,
\]
where we have used Gagliardo's trace theorem, see \cite{gagliardo}. See also \cite{necas}, Theorem 4.2 (for an English recent text see, for example, the Theorem III.2.21 in \cite{Bo-Fa}). It follows
that
\begin{equation}\label{estimate-I1}
\begin{split}
\left|\int_{\p\Omega}\bn\cdot\n\bom\cdot\bom dS\right|\leq C(\e)\|\bom\|^2_{L^2(\Omega)}+\e\,\|\n\bom\|^2_{L^2(\Omega)}\,,
\end{split}
\end{equation}
for all $0<\e<1\,.$
\vskip 0.1 true cm
\noindent\textit{Control of $I_2$:}
First, one has
\begin{equation*}
\begin{split}
\int_{\Omega}\bom\cdot\n\bu\cdot\bom dx=&\sum_{i,j,k}\int_{\hat{\Omega}}\hat{\omega}_iH^{-1}_i\p_{q_i}\left(\hat{u}_j\hat{e}_j\right)\cdot(\hat{\omega}_k\hat{e}_k)
H_1H_2H_3dq_1dq_2dq_3\\
=&\sum_{i,j,k}\int_{\hat{\Omega}}\hat{u}_j\hat{\omega}_i\hat{\omega}_k(\p_{q_i}\hat{e}_j\cdot\hat{e}_k)\,H^{-1}_iH_1H_2H_3dq_1dq_2dq_3\\
&+\sum_{i,j}\int_{\hat{\Omega}}\hat{\omega}_i\,(\p_{q_i}\hat{u}_j\,)\,\hat{\omega}_jH^{-1}_iH_1H_2H_3dq_1dq_2dq_3\\
:=&I_{21}+I_{22}\,.
\end{split}
\end{equation*}
For $I_{21}$, from \eqref{eq:assumption}, one has
\begin{equation}\label{estimate-I21}
\begin{split}
|I_{21}|\leq&C\,\int_{\Omega}|\bu|\,|\bom|^2 \,dx\\
\leq&C\,\|\bu\|_{L^2(\Omega)}\|\bom\|^2_{L^4(\Omega)}\\
\leq&C\,\|\bu\|_{L^2(\Omega)}\|\bom\|^{\f12}_{L^2(\Omega)}(\|\bom\|_{L^2(\Omega)}+\|\n\bom\|_{L^2(\Omega)})^{\f32}\\
\leq&C\,\|\bu\|_{L^2(\Omega)}\|\bom\|^{2}_{L^2(\Omega)}+C\,\|\bu\|_{L^2(\Omega)}\|\bom\|^{\f12}_{L^2(\Omega)}\|\n\bom\|_{L^2(\Omega)}^{\f32}\\
\leq&C\,\|\bu\|_{L^2(\Omega)}\|\bom\|^{2}_{L^2(\Omega)}+C(\e)\|\bu\|^4_{L^2(\Omega)}\|\bom\|^{2}_{L^2(\Omega)}+\e\,\|\n\bom\|^{2}_{L^2(\Omega)}\,.
\end{split}
\end{equation}
For $I_{22}$, we consider separately the three cases $j\neq 3$; $j=3$ and $i\neq 3$; $i=j=3$.
\vskip 0.1 true cm
\noindent\textbf{Case I: $j\neq 3$}. By integration by parts, one has
\begin{equation}\label{case1}
\begin{split}
&\int_{\hat{\Omega}}\hat{\omega}_i\,(\p_{q_i}\hat{u}_j\,)\,\hat{\omega}_jH^{-1}_iH_1H_2H_3dq_1dq_2dq_3\\
=&-\int_{\hat{\Omega}}\hat{u}_j\,(\p_{q_i}\hat{\omega}_i)\,\hat{\omega}_jH^{-1}_iH_1H_2H_3dq_1dq_2dq_3
-\int_{\hat{\Omega}}\hat{u}_j\hat{\omega}_i\,(\p_{q_i}\hat{\omega}_j)\,H^{-1}_iH_1H_2H_3dq_1dq_2dq_3\\
&-\int_{\hat{\Omega}}\hat{u}_j\hat{\omega}_i\hat{\omega}_j\p_{q_i}\left(H^{-1}_iH_1H_2H_3\right)dq_1dq_2dq_3\,.
\end{split}
\end{equation}
\noindent\textbf{Case II: $j=3$ and $i\neq 3$}.
From \eqref{eq:curl} one has
\[
\hat{\omega}_3=\f{1}{H_1H_2}\left(\f{\p(\hat{u}_2H_2)}{\p q_1}-\f{\p(\hat{u}_1H_1)}{\p q_2}\right)\,.
\]
Hence, by integration by parts, it follows that
\begin{equation}\label{case2}
\begin{split}
&\int_{\hat{\Omega}}\hat{\omega}_i\,(\p_{q_i}\hat{u}_3)\,\hat{\omega}_3H^{-1}_iH_1H_2H_3dq_1dq_2dq_3\\
=&\int_{\hat{\Omega}}
\hat{\omega}_i\,(\p_{q_i}\hat{u}_3)\,\big(\p_{q_1}(H_2\hat{u}_2)-\p_{q_2}(H_1\hat{u}_1)\big)H^{-1}_iH_3dq_1dq_2dq_3\\
=&-\int_{\hat{\Omega}}\hat{u}_2\,\p_{q_1}\big(\,\hat{\omega}_i\,(\p_{q_i}\hat{u}_3)\,H^{-1}_iH_3\,\big) H_2dq_1dq_2dq_3\\
&+\int_{\hat{\Omega}}\hat{u}_1\p_{q_2}\,\big(\hat{\omega}_i\,(\p_{q_i}\hat{u}_3)\,H^{-1}_iH_3\big)\,H_1dq_1dq_2dq_3\,.
\end{split}
\end{equation}
\noindent\textbf{Case III: $i=j=3$}. Note that, due to $\n\cdot\bu=0$, it follows
\begin{equation}\label{divz}
\f{\p(\hat{u}_1H_2H_3)}{\p q_1}+\f{\p(\hat{u}_2H_1H_3)}{\p q_2}+\f{\p(\hat{u}_3H_1H_2)}{\p q_3}=0\,.
\end{equation}
One has
\begin{equation}\label{case3}
\begin{split}
&\int_{\hat{\Omega}}\hat{\omega}_3\,(\p_{q_3}\hat{u}_3)\,\hat{\omega}_3H_1H_2dq_1dq_2dq_3\\
=&\int_{\hat{\Omega}}\hat{\omega}_3\p_{q_3}\left(H_1H_2\hat{u}_3\right)\hat{\omega}_3dq_1dq_2dq_3
-\int_{\hat{\Omega}}\hat{\omega}_3\hat{u}_3\hat{\omega}_3\p_{q_3}(H_1H_2)dq_1dq_2dq_3\\
=&-\int_{\hat{\Omega}}\hat{\omega}_3\p_{q_1}\left(H_2H_3\hat{u}_1\right)\hat{\omega}_3dq_1dq_2dq_3
-\int_{\hat{\Omega}}\hat{\omega}_3\p_{q_2}\left(H_1H_3\hat{u}_2\right)\hat{\omega}_3dq_1dq_2dq_3\\
&-\int_{\hat{\Omega}}\hat{\omega}_3\hat{u}_3\hat{\omega}_3\p_{q_3}(H_1H_2)dq_1dq_2dq_3\\
=&\int_{\hat{\Omega}}\hat{u}_1\p_{q_1}\left(\hat{\omega}^2_3\right)H_2H_3dq_1dq_2dq_3
+\int_{\hat{\Omega}}\hat{u}_1\p_{q_2}\left(\hat{\omega}^2_3\right)H_1H_3dq_1dq_2dq_3\\
&-\int_{\hat{\Omega}}\hat{u}_3\hat{\omega}^2_3\p_{q_3}(H_1H_2)dq_1dq_2dq_3\,,
\end{split}
\end{equation}
where the first equality is an identity, the second is obtained by appealing to \eqref{divz}, and the third one
follows by integration by parts.
From \eqref{case1}, \eqref{case2}, \eqref{case3} and the assumption \eqref{eq:assumption}, one can obtain
\begin{equation*}
\begin{split}
|I_{22}|\leq C\int_{\Omega}|\bar{\bu}||\n\bu||\n^2\bu|dx+C\int_{\Omega}|\bu||\n\bu|^2dx+C\int_{\Omega}|\bu|^2|\n\bu|dx\,.
\end{split}
\end{equation*}
It is easy to get that
\begin{equation*}
\begin{split}
\int_{\Omega}|\bu|^2|\n\bu|dx\leq&\|\bu\|_{L^2(\Omega)}\|\bu\|_{L^4(\Omega)}\|\n\bu\|_{L^4(\Omega)}\\
\leq&\|\bu\|^2_{L^2(\Omega)}\|\bu\|^2_{L^4(\Omega)}+\|\n\bu\|^2_{L^4(\Omega)}\\
\leq&\|\bu\|^2_{L^2(\Omega)}\|\bu\|^2_{L^4(\Omega)}+\|\n\bu\|^{\f12}_{L^2(\Omega)}(\|\n\bu\|_{L^2(\Omega)}+\|\n^2\bu\|_{L^2(\Omega)})^{\f32}\\
\leq&\|\bu\|^2_{L^2(\Omega)}\|\bu\|^2_{L^4(\Omega)}+C(\epsilon)\|\n\bu\|^2_{L^2(\Omega)}+\epsilon\|\n^2\bu\|^2_{L^2(\Omega)}\,,
\end{split}
\end{equation*}
and similarly to the proof of  \eqref{estimate-I21}
\begin{equation*}
\begin{split}
\int_{\Omega}|\bu||\n\bu|^2dx\leq C\,\|\bu\|_{L^2(\Omega)}\|\n\bu\|^{2}_{L^2(\Omega)}+C(\e)\|\bu\|^4_{L^2(\Omega)}\|\n\bu\|^{2}_{L^2(\Omega)}+\e\,\|\n^2\bu\|^{2}_{L^2(\Omega)}\,.
\end{split}
\end{equation*}
Hence, one has
\begin{equation}\label{estimate-I22}
\begin{split}
|I_{22}|\leq& C\int_{\Omega}|\bar{\bu}||\n\bu||\n^2\bu|dx+C\|\bu\|^2_{L^2(\Omega)}\|\bu\|^2_{L^4(\Omega)}+C(\e)\|\n\bu\|^2_{L^2(\Omega)}\\
&+C\,\|\bu\|_{L^2(\Omega)}\|\n\bu\|^{2}_{L^2(\Omega)}+C(\e)\|\bu\|^4_{L^2(\Omega)}\|\n\bu\|^{2}_{L^2(\Omega)}+C\e\,\|\n^2\bu\|^{2}_{L^2(\Omega)}\,.
\end{split}
\end{equation}
By H\"{o}lder's inequality, interpolation, and a Sobolev's embedding theorem, one can easily show that
\begin{equation}\label{eq:baru}
\begin{split}
&\int_{\Omega}|\bar{\bu}||\n\bu||\n^2\bu|dx\\
\leq&\||\bar{\bu}|\n\bu\|_{L^2(\Omega)}\|\n^2\bu\|_{L^2(\Omega)}\\
\leq&\|\bar{\bu}\|_{L^{p}(\Omega)}\|\n\bu\|_{L^{\f{2p}{p-2}}(\Omega)}\|\n^2\bu\|_{L^2(\Omega)}\\
\leq&\|\bar{\bu}\|_{L^{p}(\Omega)}\|\n\bu\|^{1-\f{3}{p}}_{L^{2}(\Omega)}\|\n\bu\|^{\f3{p}}_{L^6(\Omega)}\|\n^2\bu\|_{L^2(\Omega)}\\
\leq&C\|\bar{\bu}\|_{L^{p}(\Omega)}\|\n\bu\|^{1-\f{3}{p}}_{L^{2}(\Omega)}
(\|\n\bu\|_{L^2(\Omega)}+\|\n^2\bu\|_{L^2(\Omega)})^{\f3{p}}\|\n^2\bu\|_{L^2(\Omega)}\\
\leq&C\left(\|\bar{\bu}\|_{L^{p}(\Omega)}\|\n\bu\|_{L^2(\Omega)}\|\n^2\bu\|_{L^2(\Omega)}
+\|\bar{\bu}\|_{L^{p}(\Omega)}\|\n\bu\|^{1-\f{3}{p}}_{L^{2}(\Omega)}\|\n^2\bu\|_{L^2(\Omega)}^{1+\f3{p}}\right)\\
\leq&C(\e)\left(\|\bar{\bu}\|^{\f{2p}{p-3}}_{L^{p}(\Omega)}
+\|\bar{\bu}\|^2_{L^{p}(\Omega)}\right)\|\n\bu\|^{2}_{L^{2}(\Omega)}+\,\e\,\|\n^2\bu\|^{2}_{L^2(\Omega)}\,.
\end{split}
\end{equation}
Collecting \eqref{eq:Basicestimate} and the estimates \eqref{estimate-I1}, \eqref{estimate-I21}, \eqref{estimate-I22} and \eqref{eq:baru}, one obtains
\begin{equation*}
\begin{split}
\f12\p_t\int_{\Omega}|\bom|^2dx+\int_{\Omega}|\n\bom|^2dx\leq& C(\e)\left(1+\|\bu\|_{L^2(\Omega)}+\|\bu\|^4_{L^2(\Omega)}+\|\bar{\bu}\|^{\f{2p}{p-3}}_{L^{p}(\Omega)}
+\|\bar{\bu}\|^2_{L^{p}(\Omega)}\right)\|\n\bu\|^{2}_{L^{2}(\Omega)}\\
&+C\|\bu\|^2_{L^2(\Omega)}\|\bu\|^2_{L^4(\Omega)}+C\e\,\|\n^2\bu\|^{2}_{L^2(\Omega)}\,.
\end{split}
\end{equation*}
\par%
On the other hand, the following well known estimates (see for instance Theorem IV.4.8 and Theorem IV.4.9 in \cite{Bo-Fa}), hold:
\begin{equation}\label{eq:div-curl}
\|\n\bu\|_{L^2(\Omega)}\leq C\|\bom\|_{L^2(\Omega)}\,,\quad \|\n^2\bu\|_{L^2(\Omega)}\leq C\left(\|\bu\|_{L^2(\Omega)}+\|\bom\|_{H^1(\Omega)}\right)\,.
\end{equation}
Therefore, from equation \eqref{eq:div-curl}, by letting $\e$ be sufficiently small, one has
\begin{equation*}
\begin{split}
\p_t\int_{\Omega}|\bom|^2dx+\int_{\Omega}|\n\bom|^2dx\leq& C\left(1+\|\bu\|_{L^2(\Omega)}+\|\bu\|^4_{L^2(\Omega)}+\|\bar{\bu}\|^{\f{2p}{p-3}}_{L^{p}(\Omega)}
+\|\bar{\bu}\|^2_{L^{p}(\Omega)}\right)\|\bom\|^{2}_{L^{2}(\Omega)}\\
&+C\|\bu\|^2_{L^2(\Omega)}\|\bu\|^2_{L^4(\Omega)}+C\|\bu\|^2_{L^2(\Omega)}\,.
\end{split}
\end{equation*}
Finally \eqref{eq:Strong} follows by taking into account  equations \eqref{eq:RCcondition} ($q\geq\f{2p}{p-3}>2$) and \eqref{eq:div-curl}, and by appealing to a well known argument, which is based on Gronwall's inequality. Recall that weak solutions verify $\|\bu\|_{L^2(\Omega)}\in \,L^{\infty}(0,\,T)$ and $\|\bu\|^2_{L^6(\Omega)}\in L^1(0,\,T)$. Hence we have proved that $\bu$ is a strong solution.
\end{proof}
\section{On related slip boundary conditions.}\label{moregen}
In this section we present a first attempt to prove the statement of Theorem \ref{thm} with the slip boundary condition \eqref{eq:bry} replaced by the slip  boundary condition \eqref{eq:bbb} (assumed in reference \cite{Beirao18}) by means of a simple modification of our proof. This attempt fails. Hence this significant problem remains open to further investigation. This leads us to briefly show our calculations.%

\vspace{0.2cm}

Let's start by explaining our guess. As still shown in Remark \ref{re:due} condition \eqref{eq:bbb} is equivalent to
\begin{equation}\label{eq:kslip}
\bu\cdot\bn=0\,,\quad (\bom\times\bn)\cdot\tau =\,2\,\kappa_{\tau}\,\bu\cdot\tau\,,\quad \text{on}\quad \p\Omega\,.
\end{equation}
We may replace the arbitrary tangent vector $\tau$ simply by a couple of independent vectors like, for instance, the principal direction's vectors $\,\tau_1$ and $\,\tau_2\,$. In this case $\kappa_1=\,\kappa_{\tau_1}$ and $\kappa_2=\,\kappa_{\tau_2}$ are the maximum and the minimum principal curvatures.\par%
A more natural choice here is to consider the couple of tangent, orthogonal, vectors $\,\hat{e}_1 \,$ and $\,\hat{e}_2\,.$ In this case $\kappa_1$ and $\kappa_2$ are the related curvatures.  This second choice easily  leads to the couple of linear equations
\begin{equation}\label{eq:syst}
\begin{cases}
\hat{\omega}_2=\,2\,\kappa_1\,\hat{u}_1\,,\\
\hat{\omega}_1=\,-\,2\,\kappa_2\,\hat{u}_2\,.
\end{cases}
\end{equation}
Hence to replace the slip boundary condition \eqref{eq:bry} by  $\,[D(\bu)\bn]\cdot\tau=\,0\,$ means to replace assumption \eqref{eq:bry1} by \begin{equation}\label{eq:bry2}
\hat{u}_3=0\,, \quad \hat{\omega}_1=-\,2\,\kappa_2 \hat{u}_2\,,\quad \hat{\omega}_2=\,2\,\kappa_1\,\hat{u}_1\,, \quad\text{as \, $q_3=\rho_0\,,\rho_1$}\,.
\end{equation}
To prove our main statement with the boundary condition \eqref{eq:bry} replaced by the boundary condition \eqref{eq:kslip} we have to control some new boundary integrals, which no longer vanish since now $\,\hat{\omega}_1\,$ and $\,\hat{\omega}_2\,$ do not vanish. However, by \eqref{eq:syst},  $\,\hat{\omega}_1\,$ and $\,\hat{\omega}_2\,$ can be expressed in terms of the (lower order) velocity components $\,\hat{u}_2\,$ and $\,\hat{u}_1\,.$ Well known inverse trace theorems allow us to control boundary-norms of these two components by suitable internal norms. Since our P-S assumption guarantees additional regularity just for these two velocity components, one could expect that the above internal norms could be estimated in a convenient way. Unfortunately this device seems not sufficient to prove our goal. So this interesting problem remains open.%

\vspace{0.2cm}
Next we pass to showing our calculations. Let's turn back to equation \eqref{eq:comes}, by taking into account that now we can not apply to $\,\hat{w}_1=\,\hat{w}_2=\,0\,$. One has
\begin{equation*}
\begin{split}
(2l-1)\int_{\p\Omega^l}\bn\cdot\n\bom\cdot\bom dS=&(2l-1)\int_{\p\Omega^l}\bn\cdot\n\left(\f{|\bom|^2}{2}\right) dS\\
=&\int_0^1\int_0^1\left[\p_{q_3}\left(\f{|\hat{\bom}|^2}{2}\right)H_1H_2H_3^{-1}\right]\Big|_{q_3=\rho_l}dq_1dq_2\\
=&\int_0^1\int_0^1\left[\sum_j \,(\p_{q_3}\hat{w}_j)\,\hat{w}_jH_1H_2H_3^{-1}\right]\Big|_{q_3=\rho_l}dq_1dq_2\,.
\end{split}
\end{equation*}
We will drop terms which could be easily manipulated, called here "lower order terms".
Dropping lower order terms and also cancelling non significant multiplication coefficients, lead us to introduce the symbols $\,"\simeq"\,$ and $\,"\preceq"\,$, which have a clear meaning here.\par%
One has
\begin{equation}\label{eq:jjj}
(\p_{q_3}\hat{w}_j)\,\hat{w}_j\, H_1 H_2 H_3^{-1}=\,\p_{q_3}\left(H_1 H_2 \hat{w}_j \right)H_3^{-1}\hat{\omega}_j
\,-\,\p_{q_3}\left(H_1H_2\right)H_3^{-1}\hat{\omega}^2_j \preceq \p_{q_3}\left(H_1 H_2 \hat{w}_j \right)H_3^{-1}\hat{\omega}_j\,.
\end{equation}
Since $\n\cdot\bom=0$, from \eqref{eq:div} one gets
\begin{equation*}
\f{\p(\hat{\omega}_1H_2H_3)}{\p q_1}+\f{\p(\hat{\omega}_2H_1H_3)}{\p q_2}+\f{\p(\hat{\omega}_3H_1H_2)}{\p q_3}=0\,,
\end{equation*}
which gives, on $\p\hat{\Omega}\,,$
\begin{equation*}
\p_{q_3}\left(H_1H_2\hat{\omega}_3\right)H_3^{-1}\hat{\omega}_3
=\,-\,\p_{q_1}\left(H_2H_3\,\hat{\omega}_1\right)H_3^{-1}\hat{\omega}_3\,
-\,\p_{q_2}\left(H_1H_3\hat{\omega}_2\right)H_3^{-1}\hat{\omega}_3\,.
\end{equation*}
Under the new boundary conditions we can not apply to $\,\hat{\omega}_1=\,\hat{\omega}_2 =\,\p_{q_1}\,\hat{\omega}_1=\,\p_{q_2}\hat{\omega}_2=\,0\,$ on $\,\p\hat{\Omega}^l\,$ to claim the cancellation of the above right hand side. By noting that the two terms on the right hand side are symmetric, with respect to the index $1$ and $2$, we may consider just the first one.\par%
One has
$$
\p_{q_1}\left(H_2 H_3\,\hat{\omega}_1\right)H_3^{-1}\hat{\omega}_3 =\, (\p_{q_1}\,\hat{\omega}_1) \hat{\omega}_3\, H_2
\,-\,\hat{\omega}_1\, \hat{\omega}_3 \,\p_{q_1}(H_2 H_3) H_3^{-1} \simeq  (\p_{q_1}\,\hat{\omega}_1) \hat{\omega}_3\,.
$$
Note that the smooth coefficients $H_j$, as their derivatives, are not significant on our estimate below. Further, since $\p_{q_1}$ is a tangential derivative, we may apply to the second equality \eqref{eq:syst} to assume that
$\,\p_{q_1}\,\hat{\omega}_1 \simeq\,- \p_{q_1}\,\hat{u}_2\,$ on $\p\hat{\Omega}\,.$ Hence
\begin{equation}\label{eq:tres}
\int_0^1\int_0^1\left[\p_{q_1}\left(H_2 H_3\,\hat{\omega}_1\right)H_3^{-1}\hat{\omega}_3\right]|_{q_3=\rho_l}dq_1dq_2 \simeq
\int_0^1\int_0^1 \,(\p_{q_1}\hat{u}_2\,)\,\hat{\omega}_3 |_{q_3=\rho_l}dq_1dq_2\,.
\end{equation}
By appealing to Gagliardo's theorem we show that the above right hand side is bounded by
$\,C(\e) \,\|\nabla\, \bu\,\|_2 +\,\e\,\|\nabla^2\, \bu\,\|_2\,,$ which is sufficient to our purposes.%

\vspace{0.2cm}

Let's now consider in equation \eqref{eq:jjj} the terms $\, \p_{q_3} (H_1 H_2 \hat{w}_j)\,H_3^{-1} \hat{w}_j\,,$ for $j=\,1,\,2\,.$ Assume, for instance, $j=\,1\,.$ One has
$\,\p_{q_3}(H_1 H_2 \hat{w}_1) H_3^{-1} \hat{w}_1 \simeq \,(\p_{q_3} \hat{w}_1) \hat{u}_2\,.$ Hence we need to control the integral
\begin{equation}\label{eq:umm}
\int_0^1\int_0^1 \,(\p_{q_3}\hat{\omega}_1\,)\,\hat{u}_2 |_{q_3=\rho_l}dq_1dq_2\,.
\end{equation}
Roughly speaking the above integrand has the same order as that on the right hand side of \eqref{eq:tres}.
However in \eqref{eq:umm} the derivation symbol $\,\p_{q_3}\,$ appears now in the "bad position". A suitable control of the above integral turns out to be not obvious.%


\begin{thebibliography}{99}
%
\bibitem{Bae97}
\newblock Bae,~H.-O., Choe,~H.J.:
\newblock $L^{\infty}$-bound of weak solutions to Navier-Stokes equations.
\newblock In: Proceedings of the Korea-Japan Partial Differential Equations Conference (Taejon, 1996). Lecture Notes Ser. 39. Seoul Nat. Univ., Seoul, p. 13 (1997)
\bibitem{Bae07}
\newblock Bae,~H.-O., Choe,~H.J.:
\newblock {A regularity criterion for the Navier-Stokes equations.}
\newblock Commun. Partial Differ. Equations \textbf{32}, 1173--1187 (2007)
\bibitem{wolf}
\newblock Bae,~H.-O., Wolf,~J.:
\newblock {A local regularity condition involving two velocity components of Serrin-type for the Navier--Stokes equations.}
\newblock  C. R. Acad. Sci. Paris, Ser. I  \textbf{354}, 167--174 (2016)
\bibitem{batchelor}
\newblock Batchelor,~G.K:
\newblock {An Introduction to Fluid Dynamics.}
\newblock Cambridge University Press, Cambridge 1967.
\bibitem{bdvjmfm}
\newblock Beir\~ao~da~Veiga,~H.:
\newblock {On the Smoothness of a Class of Weak Solutions to the Navier--Stokes equations.}
\newblock J. Math. Fluid Mech. \textbf{2}, 315--323 (2000)
\bibitem{bdvadvances}
\newblock Beir\~ao~da~Veiga,~H.:
\newblock {Regularity for Stokes and generalized Stokes systems under nonhomogeneous slip-type boundary conditions.}
\newblock Adv. Differential Equations \textbf{9}, 1079--1114 (2004).
\bibitem{Beirao17}
\newblock Beir\~{a}o da Veiga,~H.:
\newblock {On the extension to slip boundary conditions of a Bae and Choe regularity criterion
for the Navier-Stokes equations. The half space case.}
\newblock J. Math. Anal. Appl. \textbf{453} 212--220 (2017)
\bibitem{Beirao18}
\newblock Beir\~{a}o da Veiga,~H., Bemelmans,~J., Brand,~J.:
\newblock {On a two components condition for regularity of the 3D Navier-Stokes equations under physical slip boundary conditions on non-flat boundaries.}
\newblock Mathematische Annalen  1--38 (2018)
\bibitem{beirao-crispo}
\newblock Beir\~{a}o da Veiga,~H., Crispo,~F.:
\newblock {Concerning the $W^{k,p}$-inviscid limit for 3-D flows under a slip boundary condition.}
\newblock Journal of Mathematical Fluid Mechanics \textbf{13}(1), 117--135 (2011)
\bibitem{berselli}
\newblock Berselli,~L.C.:
\newblock {A note on regularity of weak solutions
                of the Navier-Stokes equations in $\R^n.$}
\newblock Jpn. J. Math. \textbf{28}, 51--60 (2002)
\bibitem{Bo-Fa}
\newblock Boyer,~F., Fabrie,~P.:
\newblock Mathematical tools for the study of the incompressible Navier-Stokes equations and related models.
\newblock Springer Science \& Business Media (2012)
\bibitem{cao-titi}
\newblock Cao,~C., Titi,~E.S.:
\newblock {Regularity Criteria for the Three-dimensional Navier--Stokes Equations.}
\newblock  Indiana Univ. Math. J. \textbf{57}, 2643--2661 (2008)
\bibitem{chae-choe}
\newblock Chae,~D., Choe,~H.-J.:
\newblock {Regularity of Solutions to the Navier-Stokes Equation.}
\newblock  Electron. J. Differential Equations \textbf{05}, 1--7 (1999)
\bibitem{seregin}
\newblock Escauriaza,~L., Seregin,~G., \v{S}ver\'ak,~V.:
\newblock {$L_{3,\infty}$-Solutions to the Navier-Stokes Equations and Backward Uniqueness.}
\newblock Russian Math. Surveys \textbf{58}, 211--250 (2003)
\bibitem{Foias}
\newblock Foias,~C.:
\newblock {Une remarque sur l'unicit\'{e} des solutions des \'{e}quations de Navier-Stokes en dimension $n$.}
\newblock Bull. Soc. Math. Fr. \textbf{89}, 1--8 (1961)
\bibitem{gagliardo}
\newblock Gagliardo,~E.:
\newblock {Caratterizzazioni delle tracce sulla frontiera relative ad alcune classi di funzioni in $n$ variabili.}
\newblock Rend. Sem. Mat. Univ. Padova \textbf{27}, 284--305 (1957)
\bibitem{galdi2}
\newblock Galdi,~G.P.:
\newblock {An Introduction to the Navier-Stokes Initial-Boundary Value Problems.}
\newblock In: Galdi,~G.P., Heywood,~M.I., Rannacher,~R. (eds.)
              Fundamental Directions in Mathematical Fluid Mechanics.
              Advances in Mathematical Fluid Mechanics, pp. 1--70,
              Birkh\"auser, Basel (2000)
\bibitem{MR91e:35066}
\newblock Galdi,~G.P., Maremonti,~P.:
\newblock {Sulla regolarit\`a delle soluzioni deboli al
                sistema di Navier-Stokes in domini arbitrari.}
\newblock Ann. Univ. Ferrara Sez. VII Sci. Mat. \textbf{34}, 59--73 (1988)
\bibitem{giga}
\newblock Giga,~Y.:
\newblock {Solutions for Semilinear Parabolic Equations in $L^p$
                and Regularity of Weak Solutions of the Navier--Stokes System.}
\newblock J. Differential Equations \textbf{62}, 186--212 (1986)
\bibitem{cheng}
\newblock He,~C.:
\newblock {Regularity for solutions to the Navier-Stokes equations
                with one velocity component regular.}
\newblock Electron. J. Differential Equations \textbf{29}, 1--13 (2002)
\bibitem{Hopf}
\newblock Hopf, E.:
\newblock {\"{U}ber die Anfangswertaufgabe f\"{u}r die hydrodynamischen Grundgleichungen.}
\newblock  Math. Nachr. \textbf{4}, 213--231 (1951)
\bibitem{kiselev}
\newblock Kiselev,~A.A., Ladyzhenskaya,~O.A.:
\newblock {On the existence and uniqueness of the solution
                of the nonstationary problem for a viscous, incompressible fluid.}
\newblock  Izv. Akad. Nauk SSSR Ser. Mat. \textbf{21}, 655--680 (1957)
\bibitem{KIR}
\newblock Kochin,~N.E., Il'ja,~A.K., Roze,~N.V.:
\newblock Theoretical hydromechanics.
\newblock Interscience. (1964)
\bibitem{sohr-2}
\newblock Kozono,~H., Sohr,~H.:
\newblock {Regularity Criterion on Weak Solutions to the Navier-Stokes Equations.}
\newblock Adv. Differential Equations \textbf{2}, 2924--2935 (2007)
\bibitem{kukavica}
\newblock Kukavica,~I., Ziane,~M.:
\newblock {Navier-Stokes equations with regularity in one direction.}
\newblock  J. Math. Phys. \textbf{48}, 2643--2661 (2007)
\bibitem{lady-2}
\newblock Ladyzhenskaya,~O.A.:
\newblock {On uniqueness and smoothness of generalized solutions to the Navier-Stokes equations.}
\newblock Zap. Nauchn. Sem. Leningrad. Otdel. Mat. Inst. Steklov. (LOMI) \textbf{5}, 169--185 (1967)
\bibitem{ladyz}
\newblock Ladyzhenskaya,~O.A.:
\newblock {La th\'eorie math\'ematique des fluides visqueux incompressibles.}
\newblock Moscou (1961). English edition. 2nd edn. Gordon \& Breach, New York (1969)
\bibitem{Leray}
\newblock Leray,~J.:
\newblock {Sur le mouvement d'un liquide visqueux emplissant l'espace.}
\newblock Acta Math. \textbf{63}, 193--248 (1934)
\bibitem{lions}
\newblock Lions,~J.L.:
\newblock {Sur l'existence de solutions des \'equations de Navier-Stokes.}
\newblock  C. R. Acad. Sci. Paris, \textbf{248}, 2847--2849 (1959)
\bibitem{shilkin}
\newblock Mikhailov,~A.S., Shilkin,~T.N.:
\newblock {$\,L_{3,\,\infty}$-solutions to the 3D-Navier--Stokes system
                in a domain with a curved boundary.}
\newblock J. Math. Sci. (N. Y.) \textbf{143}, 2924--2935 (2007)
\bibitem{Navier}
\newblock Navier,~C.L.M.H.:
\newblock {Sur les lois de l'\'equilibre et du mouvement des corps \'elastiques.}
\newblock Mem. Acad. R. Sci.Paris, \textbf{6}, 389 (1827)
\bibitem{necas}
\newblock Ne\v{c}as,~J.:
\newblock {Les M\'ethodes Directes en Th\'eorie des  \'Equations
Elliptiques.}
\newblock Academia, Prague, 1967.
\bibitem{neustupa}
\newblock Neustupa,~J., Penel,~P.:
\newblock {Anisotropic and Geometric Criteria for Interior Regularity of Weak Solutions to the
                3D Navier--Stokes Equations.}
\newblock In: Neustupa,~J., Penel,~P. (eds.) Mathematical Fluid Mechanics.
              Advances in Mathematical Fluid Mechanics, pp. 237--265.
              Birkh\"auser, Basel (2001)
\bibitem{Prodi}
\newblock Prodi,~G.:
\newblock {Un teorema di unicit\`a per le equazioni di Navier-Stokes.}
\newblock Ann. Mat. Pura Appl. \textbf{48}, 173--182 (1959)
\bibitem{prodi-2}
\newblock Prodi,~G.:
\newblock {R\'esultats r\'ecents et probl\`emes anciens dans la th\'eorie
                des \'equations de Navier-Stokes.}
\newblock In: Les \'Equations aux D\'eriv\'ees Partielles.
              Colloques Intern. du CNRS 117, pp. 181--196, Paris (1962)
\bibitem{Serrin62}
\newblock Serrin~J.:
\newblock {On the interior regularity of weak solutions of the Navier-Stokes equations.}
\newblock Arch. Ration. Mech. Anal. \textbf{9}, 187--195 (1962)
\bibitem{Serrin63}
\newblock Serrin,~J.:
\newblock {The initial value problem for the Navier-Stokes equations.}
\newblock In: Langer, R.E. (ed.) Nonlinear Problems, 69--98. University of Wisconsin Press, Madison (1963)
\bibitem{sohr}
\newblock Sohr,~H.:
\newblock {Zur Regularit\"atstheorie der instation\"aren Gleichungen von Navier-Stokes.}
\newblock Math. Z. \textbf{184}, 359--375 (1983)
\bibitem{so-sca}
\newblock Solonnikov,~V.A., \v{S}\v{c}adilov,~V.E.:
\newblock {On a boundary value problem for a stationary system of Navier-Stokes equations.}
\newblock Proc. Steklov Inst. Math. \textbf{125}, 186--199 (1973)
\bibitem{vonwahl}
\newblock von~Wahl,~W.:
\newblock {Regularity of weak solutions of the Navier-Stokes equations.}
\newblock  Proc. Sympos. Pure Math. \textbf{45}, 497--503 (1986)
\bibitem{wang}
\newblock Wang,~L., Xin,~Z., Zhang~A.:
\newblock {Vanishing viscosity limits for 3D Navier-Stokes equations with a Navier-slip boundary condition.}
\newblock J. Math. Fluid Mech. \textbf{14}, 791--825 (2012)
\bibitem{zhang}
\newblock  Zhang,~Z., Zhong,~D., Huang,~X.:
\newblock {A refined regularity criterion for the Navier-Stokes equations
                involving one non-diagonal entry of the velocity gradient.}
\newblock  J. Math. Anal. Appl. \textbf{453}, 1145--1150 (2017)
\bibitem{zhou-pokorny}
\newblock  Zhou,~Y., Pokorn\'y,~M.:
\newblock {On the regularity of the solutions of the Navier--Stokes equations
                via one velocity component.}
\newblock  Nonlinearity \textbf{23}, 1097--1107 (2010)
\end{thebibliography}
\end{document}